\documentclass[11pt]{article}
\usepackage{amssymb,amsfonts,amsmath,latexsym, epsfig,mathrsfs,colordvi, xcolor}
\usepackage{ytableau}
\usepackage{graphicx}
\usepackage{epic}
\newtheorem{theo}{Theorem}

\newtheorem{exam}[theo]{Example}

\newtheorem{coro}[theo]{Corollary}

\makeatletter \@addtoreset{equation}{section}
\@addtoreset{theo}{}\makeatother

\def\qed{\hfill \rule{4pt}{7pt}}
\def\pf{\noindent {\it Proof.} }

\def\Ln{{\mathcal{L}_n}}
\def\Lmn{{\mathcal{L}_{m,n}}}
\def\P{{\mathcal{P}}}

\begin{document}
\setlength{\unitlength}{0.4cm}

\title{Refinements of two identities on $(n,m)$-Dyck paths}

\author{Rosena R. X. Du\footnote{Corresponding Author. Email:
rxdu@math.ecnu.edu.cn.}, Kuo Yu\\ \\ Department of Mathematics, Shanghai Key
Laboratory of PMMP \\East China Normal University,
500 Dongchuan Road \\Shanghai, 200241, P. R. China.}

\date{Jan 28, 2018}
\maketitle
\noindent {\bf Abstract:}
For integers $n, m$ with $n \geq 1$ and $0 \leq m \leq n$, an $(n,m)$-Dyck path is a lattice path in the integer lattice $\mathbb{Z} \times \mathbb{Z}$ using up steps $(0,1)$ and down steps $(1,0)$ that goes from the origin $(0,0)$ to the point $(n,n)$ and contains exactly $m$ up steps below the line $y=x$. The classical Chung-Feller theorem says that the total number of $(n,m)$-Dyck path is independent of $m$ and is equal to the $n$-th Catalan number $C_n=\frac{1}{n+1}{2n \choose n}$. For any integer $k$ with $1 \leq k \leq n$, let $p_{n,m,k}$ be the total number of $(n,m)$-Dyck paths with $k$ peaks. Ma and Yeh proved that $p_{n,m,k}$=$p_{n,n-m,n-k}$ for $0 \leq m \leq n$, and $p_{n,m,k}+p_{n,m,n-k}=p_{n,m+1,k}+p_{n,m+1,n-k}$ for $1 \leq m \leq n-2$. In this paper we give bijective proofs of these two results. Using our bijections, we also get refined enumeration results on the numbers $p_{n,m,k}$ and $p_{n,m,k}+p_{n,m,n-k}$ according to the starting and ending steps.



\section{Introduction}

Let $\Ln$ denote the set of lattice paths in the integer lattice
$\mathbb{Z} \times \mathbb{Z}$ using \emph{up} steps $U=(0,1)$ and \emph{down} steps $D=(1,0)$ that go from the origin $(0,0)$ to the point $(n,n)$. We say that $n$ is the \emph{semilength} because there are $2n$ steps. It is obvious that $|\Ln|={2n \choose n}$.

For each $L\in \Ln$, we say that $L$ is an \emph{$(n,m)$-Dyck path} if $L$ contains exactly $m$ up steps below the line $y=x$. Clearly we have $0 \leq m \leq n$. When $m=0$, $L$ never passes below the line $y=x$ and is called a \emph{Dyck path} of semilength $n$. (For this reason lattice paths in $\Ln$ are sometimes called \emph{free Dyck paths} of semilength $n$ in the literature.) A nonempty Dyck path is \emph{prime} if it touches the line $y=x$ only at the starting point and the ending point. A lattice path $L \in \mathcal{L}_n$ can be considered as a word $L_1 L_2 \cdots L_{2n}$ of $2n$ letters on the alphabet $\{U,D\}$.

Let $\Lmn$ denote the set of all $(n,m)$-Dyck paths. The classical Chung-Feller theorem \cite{Chung-Feller} says that $|\Lmn|$ is independent of $m$ and is equal to the $n$-th \emph{Catalan number} $C_n=\frac{1}{n+1}{2n \choose n}$. The proof in \cite{Chung-Feller} is based on an analytic method. Narayana proved the theorem by combinatorial methods in \cite{Narayana}. Eu, Fu and Yeh studied the theorem by using the Taylor expansions of generating functions in \cite{Eu-Liu-Yeh-AAM} and obtained a refinement of this theorem in \cite{Eu-Fu-Yeh-JCTA}. Chen revisited the theorem in \cite{ChenY} by establishing a bijection.

In \cite{Ma-Yeh-EJC} Ma and Yeh studied refinements of $(n,m)$-Dyck paths by using four parameters, namely the peak, valley, double ascent and double descent. Here a \emph{peak} (resp. \emph{valley}, \emph{double ascent}, or \emph{double descent}) is two consecutive steps $UD$ (resp. $DU$, $UU$ or $DD$). Let $\P_{n,m,k}$ (resp. $\mathcal{V}_{n,m,k}$, $\mathcal{A}_{n,m,k}$, $\mathcal{D}_{n,m,k}$) denote the set of all $(n,m)$-Dyck paths with $k$ peaks (resp. valleys, double ascents, double descents). Let $\epsilon$ be a map from the set $\{U,D\}$ to itself such that $\epsilon(U)=D$ and $\epsilon(D)=U$. Furthermore, for any path $L=L_1 L_2 \cdots L_{2n} \in \mathcal{L}_n$, we define $\phi(L)=\epsilon(L_1)\epsilon(L_2)\cdots \epsilon(L_{2n})$, $\theta(L)=L_{2n}L_{2n-1}\cdots L_2 L_1$. Clearly $\phi$ is a bijection between $\mathcal{P}_{n,m,k}$ and $\mathcal{V}_{n,n-m,k}$, and $\phi \cdot \theta$ is a bijection between $\mathcal{A}_{n,m,k}$ and $\mathcal{D}_{n,m,k}$. Hence to study the refinements on these four parameters, it is sufficient to focus on peaks and double ascents only.

Let $p_{n,m,k}=|\P_{n,m,k}|$ and $a_{n,m,k}=|\mathcal{A}_{n,m,k}|$. Ma and Yeh proved in \cite{Ma-Yeh-EJC} a Chung-Feller type theorem for Dyck paths of semilength $n$ with $k$ double ascents: the total number of $(n,m)$-Dyck paths with $k$ double ascents is independent of $m$ and is equal to the \emph{Narayana number}:
\begin{equation}\label{CF_anmk}
a_{n,m,k}=\frac{1}{k+1}{n-1 \choose k}{n \choose k}, \ \ \ 0 \leq m \leq n.
\end{equation}
Ma and Yeh gave both generating function proof and bijective proof of  the above result in \cite{Ma-Yeh-EJC}. And later in \cite{Guo-WangXX} Guo and Wang generalized this result to counting lattice paths with given number of double ascents and are dominated by a cyclically shifting piece linear boundary of varying slope.

The problem of counting $(n,m)$-Dyck paths with given number of peaks turns out to be more complicated. The number $p_{n,m,k}$ is not independent of $m$, and there is even no nice formula known for $p_{n,m,k}$. However, Ma and Yeh proved in \cite{Ma-Yeh-EJC} that the following symmetric property holds for $p_{n,m,k}$:
\begin{equation}\label{p_symmetry}
p_{n,m,k}=p_{n,n-m,n-k}.
\end{equation}
Moreover, they also found the Chung-Feller property for the sum of $p_{n,m,k}$ and $p_{n,m,n-k}$ for  $1 \leq m \leq n-2$:
\begin{equation}\label{p+_symmetry}
p_{n,m,k}+p_{n,m,n-k}=p_{n,m+1,k}+p_{n,m+1,n-k}=\frac{2(n+2)}{n(n-1)}{n \choose k-1}{n \choose k+1}.
\end{equation}

The proofs for \eqref{p_symmetry} and \eqref{p+_symmetry} given in \cite{Ma-Yeh-EJC} are based on the generating function method. We believe combinatorial proofs for these two intriguing results will be very helpful in studying the properties for $(n,m)$-Dyck paths. And in this paper we give bijective proofs and refinements of equations \eqref{p_symmetry} and \eqref{p+_symmetry}.

We denote by $\P^{UD}_{n,m,k}$ the set of paths in $\P_{n,m,k}$ that start with an up step and end with a down step, and set $p^{UD}_{n,m,k}=|\P^{UD}_{n,m,k}|$. The sets $\P^{UU}_{n,m,k}$, $\P^{DU}_{n,m,k}$, and $\P^{DD}_{n,m,k}$ and the numbers $p^{UU}_{n,m,k}$, $p^{DU}_{n,m,k}$, and $p^{DD}_{n,m,k}$ are similarly defined. Let $\P^{U}_{n,m,k}=\P^{UU}_{n,m,k}\cup \P^{UD}_{n,m,k}$.

The rest of the paper is organized as follows. In section 2 we give a bijection $\Gamma$ between $\P_{n,m,k}$ and $\P_{n,n-m,n-k}$. Observe that an $(n,m)$-Dyck path is uniquely determined by the coordinates of the peaks. The idea of the bijection is to consider the $x$-coordinate set $X$ and $y$-coordinate set $Y$ of the peaks of $P$, and then take the complementary sets of $X$ and $Y$ to be the corresponding coordinate sets of the peaks of $\Gamma(P)$. From this bijection we can get not only Equation \eqref{p_symmetry}, but also the following refined results:
$p^{UU}_{n,m,k}=p^{DD}_{n,n-m,n-k},$ and $p^{UD}_{n,m,k}=p^{DU}_{n,n-m,n-k}$. In section 3 we first define an injection from $\P^U_{n, m+1, k}$ to $\P^U_{n, m, k}$, then we give a bijection between $\P^{DU}_{n,m,k} \cup \P^{UD}_{n,m,k+1}$ and $\P^{DU}_{n,m+1,k} \cup \P^{UD}_{n, m+1, k+1}$, and therefore get the following identities:
\begin{eqnarray*}
& & p^{UU}_{n,m,k}=p^{DD}_{n,m,k}=\frac{1}{n-1}\binom{n-1}{k}\binom{n-1}{k-1};\\
& & p^{DU}_{n,m,k}+p^{UD}_{n,m,k+1}=p^{DU}_{n,m+1,k}+p^{UD}_{n,m+1,k+1}=\frac{2}{n-1}\binom{n-1}{k-1}\binom{n-1}{k+1},
\end{eqnarray*}
which imply Equation \eqref{p+_symmetry}.

\section{Bijective proof and refinements for Equation \eqref{p_symmetry}}

In this section we will first define a bijection between $\P_{n,m,k}$ and $\P_{n,n-m,n-k}$.

\begin{theo}\label{bij_pnk}
For all integers $n, m, k$ with $n \geq 1$, $1 \leq k \leq n$ and $0 \leq m \leq n$, there is a bijection $\Gamma$ between $\P_{n,m,k}$ and $\P_{n,n-m,n-k}$ such that if $P \in \P^{UU}_{n,m,k}$ (resp. $P \in \P^{UD}_{n,m,k}$) then $\Gamma(P) \in \P^{DD}_{n,n-m,n-k}$ (resp. $\Gamma(P) \in \P^{DU}_{n,n-m,n-k}$).
\end{theo}
\pf Let $P \in \P_{n,m,k}$.  Suppose the coordinates of the $k$ peaks of $P$ are $(x_i,y_i), 1 \leq i \leq k$ with $0 \leq x_1<x_2< \cdots <x_k \leq n-1$ and $1 \leq y_1<y_2< \cdots <y_k \leq n$. It is easy to check that $P$ is uniquely determined by the sets $X=\{x_1,x_2, \ldots, x_k\}$ and $Y=\{y_1, y_2, \ldots, y_k\}$. Now we set $X^{\prime}=\{x^{\prime}_1, x^{\prime}_2, \ldots, x^{\prime}_{n-k}\}_{\leq}=\{0,1, \ldots, n-1\} \backslash X$, and $Y^{\prime}=\{y^{\prime}_1, y^{\prime}_2, \ldots, y^{\prime}_{n-k}\}_{\leq}=\{1, 2, \ldots, n\} \backslash Y$.  Let $\Gamma(P)$ be the lattice path which has peaks $(x^{\prime}_i,y^{\prime}_i), 1 \leq i \leq n-k$. It is not hard to see that $\Gamma(P)$ is unique and $\Gamma(\Gamma(P))=P$. Next we will prove that $\Gamma(P) \in \P_{n,n-m,n-k}$.

We first consider the case when $m=0$. In this case $P$ is a Dyck path of semilength $n$, i.e., $P$ never goes below the line $y=x$, which implies that the coordinates of the $k$ peaks  of $P$ satisfy the following condition
\begin{equation}
\left\{
\begin{aligned}
&0=x_1<x_2< \cdots <x_k \leq n-1,\\
&1 \leq y_1<y_2< \cdots <y_k =n,\\
&x_{i+1} \leq y_{i},\ \ i=1,2,\cdots, k.
\end{aligned}
\right.
\end{equation}
Therefore for the coordinates  of peaks of $P^{\prime}$ we have
$X^{\prime},Y^{\prime} \in \{1,2, \ldots, n-1\}$ and $x^{\prime}_i \geq y^{\prime}_i$ for each $i$, $1 \leq i \leq n-k$.
Thus $P^{\prime}$ goes entirely below the line $y=x$. Hence we have
$P^{\prime} \in \P_{n,n,n-k}$.

When $m \geq 1$, $P$ intersects with the line $y=x$ at least once. Suppose $(a,a)$ is one of the intersection points, with $1 \leq a \leq n-1$. And the two peaks on both sides of $(a,a)$ along the path $P$ are $(x_i,y_i)$ and $(x_{i+1},y_{i+1})$ for some $i$, $1 \leq i \leq k-1$.  There are two cases:
\begin{itemize}
\item [1.] $x_i<a<x_{i+1}$ and $y_{i}=a$;
\item[2.] $x_{i+1}=a$ and $y_i<a<y_{i+1}$,
\end{itemize}

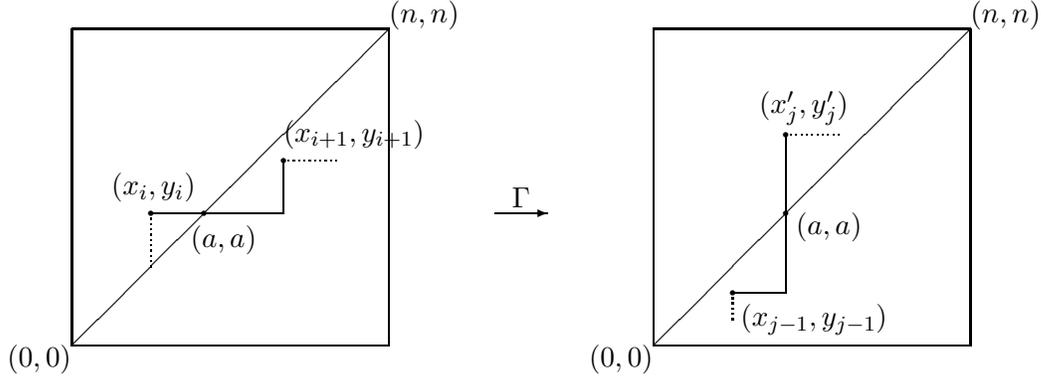
\begin{figure}

\setlength{\unitlength}{10pt}
\begin{picture}(10,10)
\put(4,0){\line(1,0){12}}
\put(4,0){\line(0,1){12}}
\put(16,0){\line(0,1){12}}
\put(4,12){\line(1,0){12}}
\put(4,0){\line(1,1){12}}
\put(7,5){\line(1,0){5}}
\bezier{10}(7,3)(7,4)(7,5)
\bezier{10}(12,7)(13,7)(14,7)
\put(12,5){\line(0,1){2}}
\put(7,5){\circle*{0.2}}
\put(9,5){\circle*{0.2}}
\put(12,7){\circle*{0.2}}
\put(5.5,6){\makebox(0,0)[l]{$(x_{i},y_{i})$}}
\put(12,8){\makebox(0,0)[l]{$(x_{i+1},y_{i+1})$}}
\put(8.5,4){\makebox(0,0)[l]{$(a,a)$}}
\put(4,-0.5){\makebox(0,0)[r]{$(0,0)$}}
\put(16,12.5){\makebox(0,0)[l]{$(n,n)$}}
\end{picture}
\put(10,5){\vector(1,0){2}}
\put(11,6){\makebox(0,0)[t]{$\Gamma$}}
\begin{picture}(10,15)
\put(16,0){\line(1,0){12}}
\put(16,0){\line(0,1){12}}
\put(28,0){\line(0,1){12}}
\put(16,12){\line(1,0){12}}
\put(16,0){\line(1,1){12}}

\bezier{5}(19,1)(19,1.5)(19,2)
\bezier{10}(21,8)(22,8)(23,8)
\put(21,2){\line(0,1){6}}
\put(19,2){\line(1,0){2}}
\put(19,2){\circle*{0.2}}
\put(21,5){\circle*{0.2}}
\put(21,8){\circle*{0.2}}
\put(19.3,1){\makebox(0,0)[l]{$(x_{j-1},y_{j-1})$}}
\put(20,9){\makebox(0,0)[l]{$(x'_{j},y'_{j})$}}
\put(21.4,4.5){\makebox(0,0)[l]{$(a,a)$}}
\put(16,-0.5){\makebox(0,0)[r]{$(0,0)$}}
\put(28,12.5){\makebox(0,0)[l]{$(n,n)$}}
\end{picture}
\begin{center}
\caption{An illustration of how $P$ and $\Gamma(P)$ intersect with the line $y=x$ for case 1.}\label{Fig_Gamma}
\end{center}
\end{figure}

In case 1, the coordinates of the peaks of $P^{\prime}$ satisfy the following condition.
\begin{equation*}
\left\{
\begin{aligned}
&0 \leq x^{\prime}_1<x^{\prime}_2< \cdots <x^{\prime}_{j}=a<x^{\prime}_{j+1} < \cdots <x^{\prime}_{n-k} \leq n-1,\\
&1 \leq y^{\prime}_1<y^{\prime}_2< \cdots <y^{\prime}_{j-1}<a<y^{\prime}_{j} < \cdots <y^{\prime}_{n-k} \leq n.\end{aligned}
\right.
\end{equation*}
with $j=a+1-i$. Hence we know the segment of path $P^{\prime}$ from $(x^{\prime}_{j-1},y^{\prime}_{j-1})$ to $(x^{\prime}_{j},y^{\prime}_{j})$ intersect with the line $y=x$ at $(a,a)$. See Figure \ref{Fig_Gamma} for an illustration.

For case 2 with similar arguments we can get the same result. Therefore we have that $P$ and $\Gamma(P)$ intersect with the line $y=x$ at the same set of points.

Now suppose $P$ intersects with the line $y=x$ at $t$ points. These points break $P$ into $t+1$ segments $P_0 P_1 \cdots P_t$. For each $i, 0 \leq i \leq t$, either $P_i \in \P_{n_i,0,k_i}$ or $P_i \in \P_{n_i,n_i,k_i}$ for some integers $n_i>0$ and $k_i\geq 0$ with $\sum_{i=0}^t n_i =n$ and $\sum_{i=0}^t k_i =k$. Moreover, from the definition of $\Gamma$ we know that $\Gamma(P)=\Gamma(P_0)\Gamma(P_1)\ldots \Gamma(P_t)$ and for each $i, 0 \leq i \leq t$ if $P_i \in \P_{n_i,0,k_i}$ then $\Gamma(P_i) \in \P_{n_i,n_i,n_i-k_i}$; if $P_i \in \P_{n_i,n_i,k_i}$ then $\Gamma(P_i) \in \P_{n_i,0,n_i-k_i}$. Therefore we proved that  $\Gamma(P) \in \P^{DD}_{n,n-m,n-k}$ (resp. $\Gamma(P) \in \P^{DU}_{n,n-m,n-k}$) if and only if $P \in \P^{UU}_{n,m,k}$ (resp. $P \in \P^{UD}_{n,m,k}$) for all integers $n, m, k$ with $n>0$, $1 \leq k \leq n$ and $0 \leq m \leq n$.
\qed

From Theorem \ref{bij_pnk} we immediately get the following result.
\begin{coro}
For all integers $n, m, k$ with $n \geq 1$, $1 \leq k \leq n$ and $0 \leq m \leq n$, we have
\begin{eqnarray}
& & p^{UU}_{n,m,k}=p^{DD}_{n,n-m,n-k};\\
& & p^{UD}_{n,m,k}=p^{DU}_{n,n-m,n-k};\\
& & p_{n,m,k}=p_{n,n-m,n-k}.
\end{eqnarray}
\end{coro}

\begin{exam}
Figure \ref{Fig_egGamma} shows an example of $P \in \P^{UD}_{10,3,5}$, and the coordinates of the 5 peaks of $P$ are $(0,2)$, $(1,4)$, $(6,6)$, $(7,9)$ and $(8,10)$. The corresponding $\Gamma(P)$ is also shown with peaks $(2,1)$, $(3,3)$, $(4,5)$, $(5,7)$ and $(9,8)$, and $\Gamma(P)\in \P^{DU}_{10,7,5}$.
\begin{figure}

\setlength{\unitlength}{15pt}
\begin{picture}(10,10)(-2,0)

\linethickness{0.1pt}
\put(0,0){\vector(1,0){11.5}}
\put(0,0){\vector(0,1){11.5}}
\dashline{0.2}(1,0)(1,10)
\dashline{0.2}(2,0)(2,10)
\dashline{0.2}(3,0)(3,10)
\dashline{0.2}(4,0)(4,10)
\dashline{0.2}(5,0)(5,10)
\dashline{0.2}(6,0)(6,10)
\dashline{0.2}(7,0)(7,10)
\dashline{0.2}(8,0)(8,10)
\dashline{0.2}(9,0)(9,10)
\dashline{0.2}(10,0)(10,10)
\dashline{0.2}(0,1)(10,1)
\dashline{0.2}(0,2)(10,2)
\dashline{0.2}(0,3)(10,3)
\dashline{0.2}(0,4)(10,4)
\dashline{0.2}(0,5)(10,5)
\dashline{0.2}(0,6)(10,6)
\dashline{0.2}(0,7)(10,7)
\dashline{0.2}(0,8)(10,8)
\dashline{0.2}(0,9)(10,9)
\dashline{0.2}(0,10)(10,10)
\thinlines
\bezier{200}(0,0)(5,5)(10,10)
\linethickness{0.9pt}
\put(0,0){\line(0,1){2}}
\put(0,2){\line(1,0){1}}
\put(1,2){\line(0,1){2}}
\put(1,4){\line(1,0){5}}
\put(6,4){\line(0,1){2}}
\put(6,6){\line(1,0){1}}
\put(7,6){\line(0,1){3}}
\put(7,9){\line(1,0){1}}
\put(8,9){\line(0,1){1}}
\put(8,10){\line(1,0){2}}
\thinlines
\put(0,2){\circle*{0.2}}
\put(1,4){\circle*{0.2}}
\put(6,6){\circle*{0.2}}
\put(7,9){\circle*{0.2}}
\put(8,10){\circle*{0.2}}
\put(0,-0.2){\makebox(0,0)[t]{$0$}}
\put(1,-0.2){\makebox(0,0)[t]{$1$}}
\put(2,-0.2){\makebox(0,0)[t]{$2$}}
\put(3,-0.2){\makebox(0,0)[t]{$3$}}
\put(4,-0.2){\makebox(0,0)[t]{$4$}}
\put(5,-0.2){\makebox(0,0)[t]{$5$}}
\put(6,-0.2){\makebox(0,0)[t]{$6$}}
\put(7,-0.2){\makebox(0,0)[t]{$7$}}
\put(8,-0.2){\makebox(0,0)[t]{$8$}}
\put(9,-0.2){\makebox(0,0)[t]{$9$}}
\put(10,-0.2){\makebox(0,0)[t]{$10$}}
\put(-0.2,1){\makebox(0,0)[r]{$1$}}
\put(-0.2,2){\makebox(0,0)[r]{$2$}}
\put(-0.2,3){\makebox(0,0)[r]{$3$}}
\put(-0.2,4){\makebox(0,0)[r]{$4$}}
\put(-0.2,5){\makebox(0,0)[r]{$5$}}
\put(-0.2,6){\makebox(0,0)[r]{$6$}}
\put(-0.2,7){\makebox(0,0)[r]{$7$}}
\put(-0.2,8){\makebox(0,0)[r]{$8$}}
\put(-0.2,9){\makebox(0,0)[r]{$9$}}
\put(-0.2,10){\makebox(0,0)[r]{$10$}}
\end{picture}
\put(4,5){\vector(1,0){2}}
\put(5,5.8){\makebox(0,0)[t]{$\Gamma$}}
\begin{picture}(10,10)(-8,0)
\linethickness{0.1pt}
\put(0,0){\vector(1,0){11.5}}
\put(0,0){\vector(0,1){11.5}}
\dashline{0.2}(1,0)(1,10)
\dashline{0.2}(2,0)(2,10)
\dashline{0.2}(3,0)(3,10)
\dashline{0.2}(4,0)(4,10)
\dashline{0.2}(5,0)(5,10)
\dashline{0.2}(6,0)(6,10)
\dashline{0.2}(7,0)(7,10)
\dashline{0.2}(8,0)(8,10)
\dashline{0.2}(9,0)(9,10)
\dashline{0.2}(10,0)(10,10)
\dashline{0.2}(0,1)(10,1)
\dashline{0.2}(0,2)(10,2)
\dashline{0.2}(0,3)(10,3)
\dashline{0.2}(0,4)(10,4)
\dashline{0.2}(0,5)(10,5)
\dashline{0.2}(0,6)(10,6)
\dashline{0.2}(0,7)(10,7)
\dashline{0.2}(0,8)(10,8)
\dashline{0.2}(0,9)(10,9)
\dashline{0.2}(0,10)(10,10)
\thinlines
\bezier{200}(0,0)(5,5)(10,10)
\linethickness{0.9pt}
\put(0,0){\line(1,0){2}}
\put(2,0){\line(0,1){1}}
\put(2,1){\line(1,0){1}}
\put(3,1){\line(0,1){2}}
\put(3,3){\line(1,0){1}}
\put(4,3){\line(0,1){2}}
\put(4,5){\line(1,0){1}}
\put(5,5){\line(0,1){2}}
\put(5,7){\line(1,0){4}}
\put(9,7){\line(0,1){1}}
\put(9,8){\line(1,0){1}}
\put(10,8){\line(0,1){2}}
\thinlines
\put(2,1){\circle*{0.2}}
\put(3,3){\circle*{0.2}}
\put(4,5){\circle*{0.2}}
\put(5,7){\circle*{0.2}}
\put(9,8){\circle*{0.2}}
\put(0,-0.2){\makebox(0,0)[t]{$0$}}
\put(1,-0.2){\makebox(0,0)[t]{$1$}}
\put(2,-0.2){\makebox(0,0)[t]{$2$}}
\put(3,-0.2){\makebox(0,0)[t]{$3$}}
\put(4,-0.2){\makebox(0,0)[t]{$4$}}
\put(5,-0.2){\makebox(0,0)[t]{$5$}}
\put(6,-0.2){\makebox(0,0)[t]{$6$}}
\put(7,-0.2){\makebox(0,0)[t]{$7$}}
\put(8,-0.2){\makebox(0,0)[t]{$8$}}
\put(9,-0.2){\makebox(0,0)[t]{$9$}}
\put(10,-0.2){\makebox(0,0)[t]{$10$}}
\put(-0.2,1){\makebox(0,0)[r]{$1$}}
\put(-0.2,2){\makebox(0,0)[r]{$2$}}
\put(-0.2,3){\makebox(0,0)[r]{$3$}}
\put(-0.2,4){\makebox(0,0)[r]{$4$}}
\put(-0.2,5){\makebox(0,0)[r]{$5$}}
\put(-0.2,6){\makebox(0,0)[r]{$6$}}
\put(-0.2,7){\makebox(0,0)[r]{$7$}}
\put(-0.2,8){\makebox(0,0)[r]{$8$}}
\put(-0.2,9){\makebox(0,0)[r]{$9$}}
\put(-0.2,10){\makebox(0,0)[r]{$10$}}
\end{picture}
\begin{center}
\caption{An example of the bijection $\Gamma$ with $P \in \P_{10,3,5}$ and $\Gamma(P) \in \P_{10,7,5}$.}\label{Fig_egGamma}
\end{center}
\end{figure}
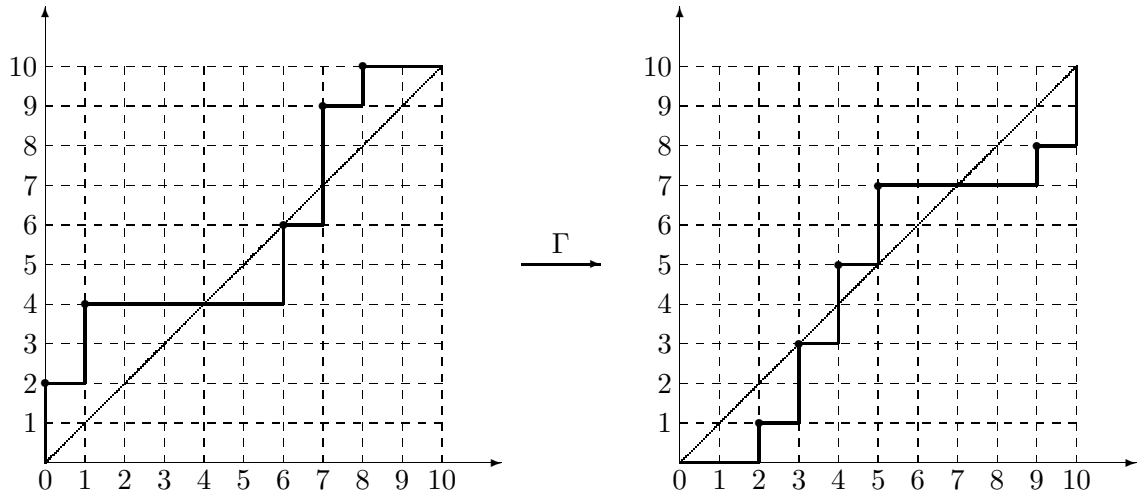
\end{exam}

\section{Bijective proof and refinements for Equation \eqref{p+_symmetry}}

In this section we first define a map $f$ from $\P^U_{n,m+1,k}$ to $\P^U_{n,m,k}$ for each $m$, $1 \leq m \leq n-2$.

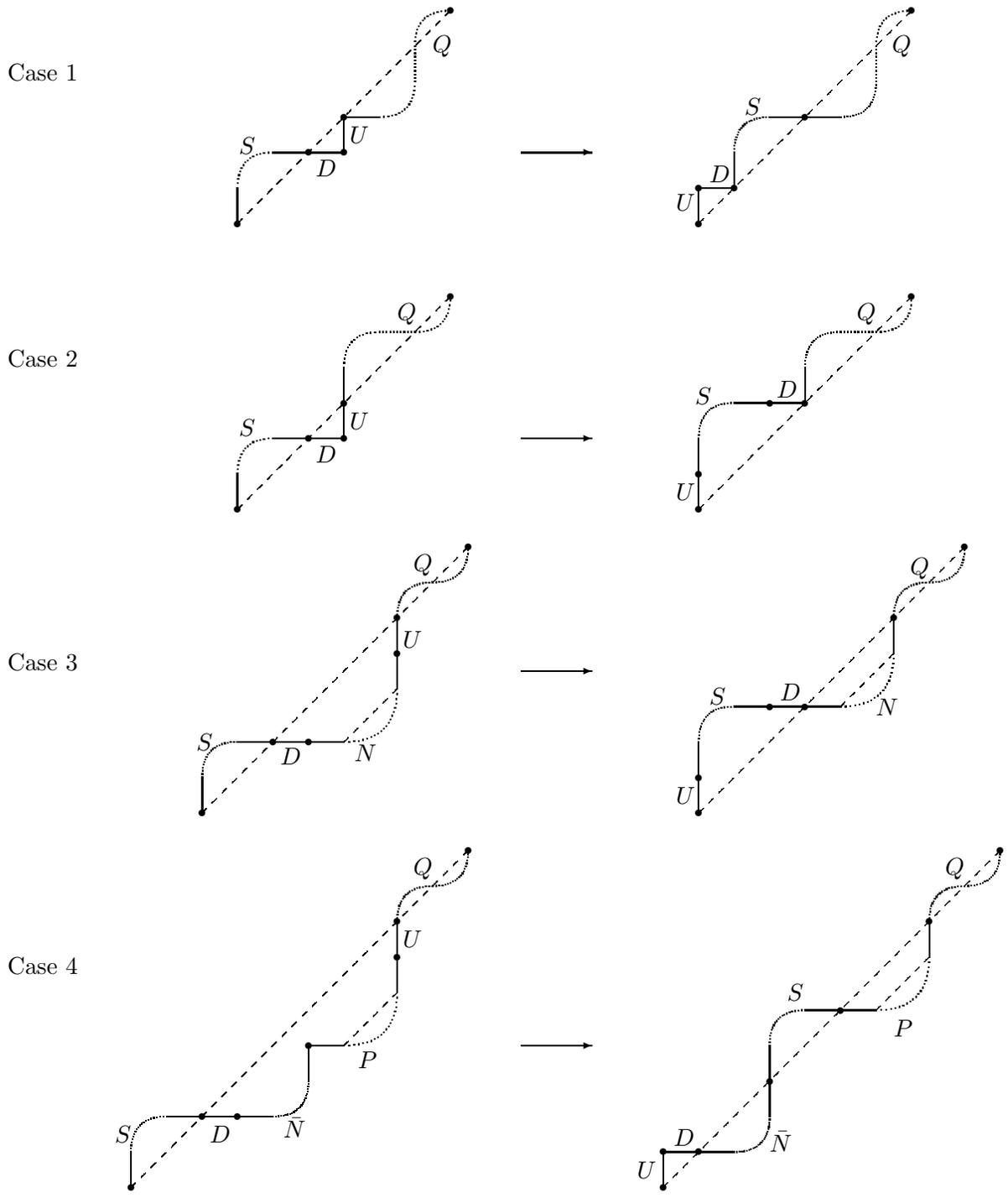
\begin{figure}
\setlength{\unitlength}{16pt}

\begin{picture}(10,6)(-3,0)
\put(-1,5){\makebox(0.5,0.5)[tr]{Case 1}}
\put(4,1){\line(0,1){1}}
\bezier{15}(4,2)(4,3)(5,3)
\put(5,3){\line(1,0){2}}
\put(7,3){\line(0,1){1}}
\put(7,4){\line(1,0){1}}
\bezier{15}(8,4)(9,4)(9,5)
\dottedline{0.1}(9,5)(9,6)
\bezier{15}(9,6)(9,7)(10,7)
\put(4,1){\circle*{0.2}}
\put(6,3){\circle*{0.2}}
\put(7,3){\circle*{0.2}}
\put(7,4){\circle*{0.2}}
\put(10,7){\circle*{0.2}}
\put(4.5,3.2){\makebox(0,0)[r]{$S$}}
\put(6.5,2.8){\makebox(0,0)[t]{$D$}}
\put(7.7,3.5){\makebox(0,0)[r]{$U$}}
\put(9.5,6){\makebox(0,0)[l]{$Q$}}
\linethickness{0.1pt}
\dashline{0.2}(4,1)(10,7)
\thinlines
\put(12,3){\vector(1,0){2}}
\put(17,1){\line(0,1){1}}
\put(17,2){\line(1,0){1}}
\put(18,2){\line(0,1){1}}
\bezier{20}(18,3)(18,4)(19,4)
\put(19,4){\line(1,0){2}}
\bezier{20}(21,4)(22,4)(22,5)
\dottedline{0.1}(22,5)(22,6)
\bezier{15}(22,6)(22,7)(23,7)
\put(17,1){\circle*{0.2}}
\put(17,2){\circle*{0.2}}
\put(18,2){\circle*{0.2}}
\put(20,4){\circle*{0.2}}
\put(23,7){\circle*{0.2}}
\put(16.9,1.6){\makebox(0,0)[r]{$U$}}
\put(17.3,2.4){\makebox(0,0)[l]{$D$}}
\put(18.3,4.3){\makebox(0,0)[l]{$S$}}
\put(23,6){\makebox(0,0)[r]{$Q$}}
\linethickness{0.1pt}
\dashline{0.2}(17,1)(23,7)
\thinlines
\end{picture}

\begin{picture}(10,8)(-3,0)
\put(-1,5){\makebox(0.5,0.5)[tr]{Case 2}}
\put(4,1){\line(0,1){1}}
\bezier{15}(4,2)(4,3)(5,3)
\put(5,3){\line(1,0){2}}
\put(7,3){\line(0,1){1}}
\put(7,4){\line(0,1){1}}
\bezier{15}(7,5)(7,6)(8,6)
\dottedline{0.1}(8,6)(9,6)
\bezier{15}(9,6)(10,6)(10,7)
\put(4,1){\circle*{0.2}}
\put(6,3){\circle*{0.2}}
\put(7,3){\circle*{0.2}}
\put(7,4){\circle*{0.2}}
\put(10,7){\circle*{0.2}}
\put(4.5,3.3){\makebox(0,0)[r]{$S$}}
\put(6.5,2.8){\makebox(0,0)[t]{$D$}}
\put(7.7,3.5){\makebox(0,0)[r]{$U$}}
\put(8.5,6.5){\makebox(0,0)[l]{$Q$}}
\linethickness{0.1pt}
\dashline{0.2}(4,1)(10,7)
\thinlines
\put(12,3){\vector(1,0){2}}
\put(17,1){\line(0,1){2}}
\bezier{20}(17,3)(17,4)(18,4)
\put(18,4){\line(1,0){2}}
\put(20,4){\line(0,1){1}}
\bezier{20}(20,5)(20,6)(21,6)
\dottedline{0.1}(21,6)(22,6)
\bezier{15}(22,6)(23,6)(23,7)
\put(17,1){\circle*{0.2}}
\put(17,2){\circle*{0.2}}
\put(19,4){\circle*{0.2}}
\put(20,4){\circle*{0.2}}
\put(23,7){\circle*{0.2}}
\put(16.9,1.5){\makebox(0,0)[r]{$U$}}
\put(16.9,4.2){\makebox(0,0)[l]{$S$}}
\put(19.2,4.4){\makebox(0,0)[l]{$D$}}
\put(22,6.5){\makebox(0,0)[r]{$Q$}}
\linethickness{0.1pt}
\dashline{0.2}(17,1)(23,7)
\thinlines
\end{picture}

\begin{picture}(10,8.5)(-2,0)
\put(0,5){\makebox(0.5,0.5)[tr]{Case 3}}
\put(4,1){\line(0,1){1}}
\bezier{20}(4,2)(4,3)(5,3)
\put(5,3){\line(1,0){3}}
\bezier{20}(8,3)(9.5,3)(9.5,4.5)

\put(9.5,4.5){\line(0,1){2}}
\bezier{20}(9.5,6.5)(9.5,7.5)(10.5,7.5)
\bezier{20}(10.5,7.5)(11.5,7.5)(11.5,8.5)
\put(4,1){\circle*{0.2}}
\put(6,3){\circle*{0.2}}
\put(7,3){\circle*{0.2}}
\put(9.5,5.5){\circle*{0.2}}
\put(9.5,6.5){\circle*{0.2}}
\put(11.5,8.5){\circle*{0.2}}
\put(4.3,3){\makebox(0,0)[r]{$S$}}
\put(6.8,2.6){\makebox(0,0)[r]{$D$}}
\put(8.3,2.7){\makebox(0,0)[l]{$N$}}
\put(10.2,5.9){\makebox(0,0)[r]{$U$}}
\put(10.5,8){\makebox(0,0)[r]{$Q$}}
\linethickness{0.1pt}
\dashline{0.2}(4,1)(11.5,8.5)
\dashline{0.2}(8,3)(9.5,4.5)
\thinlines
\put(13,5){\vector(1,0){2}}
\put(18,1){\line(0,1){2}}
\bezier{20}(18,3)(18,4)(19,4)
\put(19,4){\line(1,0){3}}
\bezier{20}(22,4)(23.5,4)(23.5,5.5)
\put(23.5,5.5){\line(0,1){1}}
\bezier{20}(23.5,6.5)(23.5,7.5)(24.5,7.5)
\bezier{20}(24.5,7.5)(25.5,7.5)(25.5,8.5)
\put(18,1){\circle*{0.2}}
\put(18,2){\circle*{0.2}}
\put(20,4){\circle*{0.2}}
\put(21,4){\circle*{0.2}}
\put(23.5,6.5){\circle*{0.2}}
\put(25.5,8.5){\circle*{0.2}}
\put(17.9,1.5){\makebox(0,0)[r]{$U$}}
\put(18.3,4.3){\makebox(0,0)[l]{$S$}}
\put(20.3,4.4){\makebox(0,0)[l]{$D$}}
\put(23,4){\makebox(0,0)[l]{$N$}}
\put(24.5,7.9){\makebox(0,0)[r]{$Q$}}
\linethickness{0.1pt}
\dashline{0.2}(18,1)(25.5,8.5)
\dashline{0.2}(22,4)(23.5,5.5)
\thinlines
\end{picture}

\begin{picture}(10,10.5)
\put(2,7){\makebox(0.5,0.5)[tr]{Case 4}}
\put(4,1){\line(0,1){1}}
\bezier{20}(4,2)(4,3)(5,3)
\put(5,3){\line(1,0){3}}
\bezier{25}(8,3)(9,3)(9,4)
\put(9,4){\line(0,1){1}}
\put(9,5){\line(1,0){1}}
\bezier{20}(10,5)(11.5,5)(11.5,6.5)
\put(11.5,6.5){\line(0,1){2}}
\bezier{20}(11.5,8.5)(11.5,9.5)(12.5,9.5)
\bezier{20}(12.5,9.5)(13.5,9.5)(13.5,10.5)
\put(4,1){\circle*{0.2}}
\put(6,3){\circle*{0.2}}
\put(7,3){\circle*{0.2}}
\put(9,5){\circle*{0.2}}
\put(11.5,7.5){\circle*{0.2}}
\put(11.5,8.5){\circle*{0.2}}
\put(13.5,10.5){\circle*{0.2}}
\put(4,2.5){\makebox(0,0)[r]{$S$}}
\put(6.8,2.5){\makebox(0,0)[r]{$D$}}
\put(8.3,2.7){\makebox(0,0)[l]{$\bar{N}$}}
\put(10.4,4.6){\makebox(0,0)[l]{$P$}}
\put(12.2,8){\makebox(0,0)[r]{$U$}}
\put(12.5,10){\makebox(0,0)[r]{$Q$}}
\linethickness{0.1pt}
\dashline{0.2}(4,1)(13.5,10.5)
\dashline{0.2}(10,5)(11.5,6.5)
\thinlines
\put(15,5){\vector(1,0){2}}
\put(19,1){\line(0,1){1}}
\put(19,2){\line(1,0){2}}
\bezier{25}(21,2)(22,2)(22,3)
\put(22,3){\line(0,1){2}}
\bezier{20}(22,5)(22,6)(23,6)
\put(23,6){\line(1,0){2}}
\bezier{20}(25,6)(26.5,6)(26.5,7.5)
\put(26.5,7.5){\line(0,1){1}}
\bezier{20}(26.5,8.5)(26.5,9.5)(27.5,9.5)
\bezier{20}(27.5,9.5)(28.5,9.5)(28.5,10.5)
\put(19,1){\circle*{0.2}}
\put(19,2){\circle*{0.2}}
\put(20,2){\circle*{0.2}}
\put(22,4){\circle*{0.2}}
\put(24,6){\circle*{0.2}}
\put(26.5,8.5){\circle*{0.2}}
\put(28.5,10.5){\circle*{0.2}}
\put(18.8,1.5){\makebox(0,0)[r]{$U$}}
\put(22.5,6.5){\makebox(0,0)[l]{$S$}}
\put(19.3,2.4){\makebox(0,0)[l]{$D$}}
\put(22,2.3){\makebox(0,0)[l]{$\bar{N}$}}
\put(25.5,5.5){\makebox(0,0)[l]{$P$}}
\put(27.5,10){\makebox(0,0)[r]{$Q$}}
\linethickness{0.1pt}
\dashline{0.2}(19,1)(28.5,10.5)
\dashline{0.2}(25,6)(26.5,7.5)
\thinlines
\end{picture}

\caption{The four cases of $f: \P^U_{n,m+1,k} \mapsto \P^U_{n,m,k}$.}
\label{Fig_f}
\end{figure}

Given any path $P \in \P^U_{n,m+1,k}$, we can uniquely decompose $P$ into $P=SDNUQ$ such that: $D$ is the first down step below the line $y=x$, and $U$ is the first up step that ends at the line $y=x$. It is easy to check that $S$ and $\phi(N)$ are both Dyck paths. ($\phi$ is defined in the fourth paragraph in the introduction.) There are four cases (See Fig \ref{Fig_f}):

\begin{itemize}
\item [1.] If $N$ is empty, and $Q$ starts with a down step, then we set $f(P)=UDSQ$;

\item [2.] If $N$ is empty, and $Q$ starts with an up step, then we set $f(P)=USDQ$;

\item [3.] If $N$ is not empty, and $\phi(N)$ is a prime Dyck path, then we set $f(P)=USDNQ$; 

\item[4.] If $N$ is not empty, and $\phi(N)$ is not a prime Dyck path.  In this case we decompose $N$ into $N=\bar{N}P$ such that $\phi(P)$ is the rightmost prime Dyck path in $\phi(N)$, and set $f(P)=UD\bar{N}SPQ$. (Note that in this case $\bar{N}$ is not empty.)
\end{itemize}

It is easy to check that for each of the above four cases, $f(P) \in \P^U_{n,m,k}$. Next we will prove that $f$ is an injection by showing that each $P^{\prime} \in \P^U_{n,m,k}$ belongs to exactly one of the following five cases, and for the first four cases, there is a unique path $P\in \P^U_{n, m+1, k}$ such that $P^{\prime}=f(P)$.

For each $P^{\prime} \in \P^U_{n, m, k}$, we can uniquely decompose $P^{\prime}$ into $P^{\prime}=USDQ$ such that $D$ is the first down step that ends at the line $y=x$ (hence $S$ is a Dyck path). There are five cases. For the first four cases we define a map $f^\prime$ which gives the inverse of the four cases of $f$, and the correspondence between the cases of $f$ and $f^\prime$ is $1 \leftrightarrow c$, $2 \leftrightarrow a$,  $3 \leftrightarrow b$, $4 \leftrightarrow d$.

\begin{itemize}
\item[a.] If $S$ is not empty, and $Q$ starts with an up step, then we set $f^{\prime}(P^\prime)=SDUQ$;

\item [b.] If $S$ is not empty, and $Q$ starts with a down step, then we decompose $Q$ into $PQ'$ such that $\phi(P)$ is the leftmost prime Dyck path in $\phi(Q)$, and set $f^{\prime}(P^\prime)=SDPUQ'$;

\item [c.] If $S$ is empty, and $Q$ starts with an up step. We decompose $Q$ into $MQ'$ such that $M$ is a maximum Dyck path, which implies that $Q'$ starts with a down step. We set $f^{\prime}(P^\prime)=MDUQ'$;

\item [d.] If $S$ is empty, $Q$ is nonempty and starts with a down step, and we can decompose $Q$ into $NMPQ'$ such that $\phi(N)$ is a Dyck path, $M$ is a maximum Dyck path, $PQ^{\prime}$ is nonempty and $\phi(P)$ is a prime Dyck path, then we set $f^{\prime}(P^\prime)=MDNPUQ'$.

\item [e.] If $P^{\prime}$ is of the form $UDNM$, where $\phi(N)$ and $M$ are both Dyck paths. 
\end{itemize}

It is easy to check that if $P^{\prime}$ belongs to Case a-d, then $f(f^{\prime}(P^\prime))=P^{\prime}$, hence $f: \P^U_{n,m+1,k} \mapsto \P^U_{n,m,k}$ is an injection. The paths in $P^{\prime}\in \P^U_{n,m,k}$ that belong to Case e form a subset of $\P^{UD}_{n,m,k}$, we denote this subset as $\hat{\P}^{UD}_{n,m,k}$, and set $\tilde{\P}^{UD}_{n,m,k}=\P^{UD}_{n,m,k} \setminus \hat{\P}^{UD}_{n,m,k}$.  Therefore we have the following theorem.

\begin{theo}\label{Bij_f}
For all positive integers $n,k$ and $m$ with $1 \leq m \leq n-2$ and $1 \leq k \leq n-1$, we have
\begin{itemize}
\item[1.] $f$ is an injection from $\P^U_{n, m+1, k}$ to $\P^U_{n, m, k}$;
\item[2.] $f$ is a bijection between $\P^{UU}_{n,m+1,k}$ and $\P^{UU}_{n,m,k}$;
\item [3.] $f$ is a bijection between $\P^{UD}_{n,m+1,k}$ and $\tilde{\P}^{UD}_{n,m,k}$.
\end{itemize}
\end{theo}

For each $P \in \hat{\P}^{UD}_{n,m,k}$, we define a map $g$ from $ \hat{\P}^{UD}_{n,m,k}$ to $\P^{DU}_{n,m+1,k-1}$ as the following.

Given any path $P \in \hat{\P}^{UD}_{n,m,k}$, we can uniquely decompose $P$ into $P=UDNM$ such that $M$ and $\phi(N)$ are both Dyck paths. We set $g(P)=NMDU$.
It is clear that $g(P)\in\P^{DU}_{n,m+1,k-1}$. We denote the set $g(\hat{\P}^{UD}_{n,m,k})$ as $\hat{\P}^{DU}_{n,m+1,k-1}$, and set $\tilde{\P}^{DU}_{n,m+1,k-1}=\P^{DU}_{n,m+1,k-1} \setminus \hat{\P}^{DU}_{n,m+1,k-1}$. So $g$ is a bijection between $\hat{\P}^{UD}_{n,m,k}$ and $\hat{\P}^{DU}_{n,m+1,k-1}$.

\begin{theo} \label{Bij_Phi}
For all positive integers $n,k$ and $m$ with $1 \leq m \leq n-2$ and $1 \leq k \leq n-1$, there is a bijection $\Phi$ between the set $\P^{DU}_{n,m+1,k} \cup \P^{UD}_{n,m+1,k+1}$ and the set $\P^{DU}_{n,m,k} \cup \P^{UD}_{n,m,k+1}$.
\end{theo}

\pf From the definition of the bijection $f$ we already know that $\P^{UD}_{n,m,k+1}$ can be decomposed into $\hat{\P}^{UD}_{n,m,k+1} \cup \tilde{\P}^{UD}_{n,m,k+1}$, and $\P^{DU}_{n,m+1,k}$ can be decomposed into $\hat{\P}^{DU}_{n,m+1,k} \cup \tilde{\P}^{DU}_{n,m+1,k}$. Now we define $\Phi$: $\tilde{\P}^{DU}_{n,m+1,k} \cup \hat{\P}^{DU}_{n,m+1,k}\cup \P^{UD}_{n,m+1,k+1} \mapsto \P^{DU}_{n,m,k} \cup \hat{\P}^{UD}_{n,m,k+1}\cup\tilde{\P}^{UD}_{n,m,k+1}$ according to the following three cases:

\begin{itemize}
\item [1.] If $P\in\P^{UD}_{n,m+1,k+1}$, then $\Phi(P)=f(P)\in\tilde{\P}^{UD}_{n,m,k+1}$;
\item [2.] If $P\in\hat{\P}^{DU}_{n,m+1,k}$, then $\Phi(P)=g^{-1}(P)\in\hat{\P}^{UD}_{n,m,k+1}$;
\item [3.] If $P\in\tilde{\P}^{DU}_{n,m+1,k}$,  we set $\Phi(P)=\theta(f^{-1}(\theta(P)))$. Here $\theta(P)$ is the reverse of $P$, as defined in the introduction. Hence we have $\theta(P)\in \tilde{\P}^{UD}_{n, n-m-1, k+1}$, and  $f^{-1}(\theta(P))\in \P^{UD}_{n, n-m, k+1}$, and therefore  $\Phi(P)=\theta(f^{-1}(\theta(P)))\in \P^{DU}_{n,m,k}$.
\end{itemize}

It is easy to check that $\Phi$ is a bijection between $\P^{DU}_{n,m+1,k} \cup \P^{UD}_{n,m+1,k+1}$ and $\P^{DU}_{n,m,k} \cup \P^{UD}_{n,m,k+1}$.
\qed

The proceeding bijections lead to several refinements on the numbers $p_{n,m,k}$ as well as Equations \eqref{p_symmetry} and \eqref{p+_symmetry}.

\begin{coro} For all positive integers $n, m, k$ with $1 \leq m \leq n-1$ and $1 \leq k \leq n-1$, we have
\begin{equation}\label{pUU}
p^{UU}_{n,m,k}=p^{DD}_{n,m,k}=\frac{1}{n-1}\binom{n-1}{k}\binom{n-1}{k-1}.
\end{equation}
\end{coro}
\pf
For each $P\in \P^{UU}_{n,m,k}$,  suppose the coordinates of the peaks of $P$ are $(x_i,y_i), 1 \leq i \leq k$ with $0 \leq x_1<x_2< \cdots <x_k \leq n-1$ and $1 \leq y_1<y_2< \cdots <y_k \leq n-1$. Since $P$ both starts and ends with an up step, we have that $x_1=0$. And $P$ is uniquely determined once the peaks are given. There are ${n-1 \choose k-1}$ ways to choose the numbers $x_2, x_3, \ldots, x_k$ and ${n-1 \choose k}$ ways to choose the numbers $y_1, y_2, \ldots, y_k$. Therefore we have
\[\sum_{m=1}^{n-1}p^{UU}_{n,m,k}={n-1 \choose k}{n-1 \choose k-1}.\]
Moreover, from Theorem \ref{Bij_f} we know that $p^{UU}_{n,m,k}=p^{UU}_{n,m+1,k}$ for all $m$, $1 \leq m \leq n-2$. Therefore we have that $p^{UU}_{n,m,k}=\frac{1}{n-1}\binom{n-1}{k}\binom{n-1}{k-1}$. The case for $p^{DD}_{n,m,k}$ can be similarly proved. \qed

\begin{coro} For all positive integers $n,k$ and $m$ with $1 \leq m \leq n-2$ and $1 \leq k \leq n-1$, we have
\begin{equation}\label{2+2}
p^{DU}_{n,m,k}+p^{UD}_{n,m,k+1}=p^{DU}_{n,m+1,k}+p^{UD}_{n,m+1,k+1}=\frac{2}{n-1}\binom{n-1}{k-1}\binom{n-1}{k+1}.
\end{equation}
\end{coro}

\pf From Theorem \ref{Bij_Phi} we know that for each $m, 1 \leq m \leq n-2$ the following holds.
\[p^{DU}_{n,m,k}+p^{UD}_{n,m,k+1}=p^{DU}_{n,m+1,k}+p^{UD}_{n,m+1,k+1}=p^{DU}_{n,1,k}+p^{UD}_{n,1,k+1}.\]
It is obvious that $p^{DU}_{n,1,k}=0$, so it remains to count $p^{UD}_{n,1,k+1}$. For any $P\in\P^{UD}_{n,1,k+1}$, we can
decompose $P$ into $P=SDUQ$ such that: $S$ and $Q$ are both Dyck paths. Since $P$ starts with an up step and ends with a down step, we know that neither $S$ nor $Q$ is empty. We define $\tau$ : $\tau(P)=USDQ$, so $\tau(P)\in\P^{UD}_{n,0,k+1}$.

On the other hand, for each $P^{\prime} \in \P^{UD}_{n,0,k+1}$, we can uniquely decompose it into $P^{\prime}=USDQ$ such that $D$ is the first down step that ends at the line $y=x$ (hence $S$ and $Q$ are both Dyck paths). And we have that $\tau^{-1}(P^{\prime})=SDUQ\in \P^{UD}_{n,1,k+1}$ if and only if neither $P$ nor $Q$ is empty. Among all the paths in $\P_{n,0,k+1}$, there are $p_{n-1,0,k}$ of them with $S$ empty in the
decomposition, and $p_{n-1,0,k+1}$ of them with $S$ not empty but $Q$ is empty. Therefore we have
\[p^{UD}_{n,1,k+1}=p_{n,0,k+1}-p_{n-1,0,k}-p_{n-1,0,k+1}.\]
With the standard result (See for example, \cite{EC2}) that $p_{n,0,k}=\frac{1}{n}{n \choose k}{n\choose k-1}$, which is the Narayana number, we can get the desired result by simple computation. \qed

\begin{coro} \cite{Ma-Yeh-EJC} For all positive integers $n,k$ and $m$ with $1 \leq m \leq n-2$ and $1 \leq k \leq n-1$, we have
\begin{equation}
p_{n,m,k}+p_{n,n-m,k}=p_{n,m+1,k}+p_{n,n-m-1,k}=\frac{2(n+2)}{n(n-1)}\binom{n}{k-1}\binom{n}{k+1}.
\end{equation}
\end{coro}

\pf It is obvious that $p_{n,m,k}=p^{UU}_{n,m,k}+p^{UD}_{n,m,k}+p^{DU}_{n,m,k}+p^{DD}_{n,m,k}$.


Replace $k$ with $k-1$ in Equation \eqref{2+2}, we get a new equation. By adding the new equation to \eqref{2+2}, we get
\[p^{UD}_{n,m,k}+p^{DU}_{n,m,k}+p^{UD}_{n,m,k+1}+p^{DU}_{n,m,k-1}=p^{UD}_{n,m+1,k}+p^{DU}_{n,m+1,k}+p^{UD}_{n,m+1,k+1}+p^{DU}_{n,m+1,k-1}.\]

Since the map $\theta$ is a bijection between $\P^{DU}_{n, m, k}$ and $\P^{UD}_{n, n-m, k+1}$, we have $p^{DU}_{n, m, k}=p^{UD}_{n, n-m, k+1}$. Therefore we have
\[p^{UD}_{n,m,k}+p^{DU}_{n,m,k}+p^{DU}_{n,n-m,k}+p^{UD}_{n,n-m,k}
=p^{UD}_{n,m+1,k}+p^{DU}_{n,m+1,k}+p^{DU}_{n,n-m-1,k}+p^{UD}_{n,n-m-1,k}.\]

Combining \eqref{pUU}, \eqref{2+2} and the above equations we get
\begin{eqnarray*}
& &p_{n,m,k}+p_{n,n-m,k}=p_{n,m+1,k}+p_{n,n-m-1,k}\\
&=&\frac{2}{n-1}\binom{n-1}{k-1}\binom{n-1}{k+1}+\frac{2}{n-1}\binom{n-1}{k-2}\binom{n-1}{k}+\frac{4}{n-1}{n-1 \choose k}{n-1 \choose k-1}\\
&=&\frac{2(n+2)}{n(n-1)}\binom{n}{k-1}\binom{n}{k+1}.
\end{eqnarray*}
\qed

\vskip 2mm \noindent{\bf Acknowledgments.} This work is partially supported by the Science and Technology Commission of Shanghai Municipality (STCSM), No. 13dz2260400.


\begin{thebibliography}{99}

\bibitem{ChenY}
Y.M. Chen, The Chung-Feller theorem revisited, Discrete Math. 308 (2008) 1328--1329.

\bibitem {Chung-Feller}
K.L.Chung, W. Feller, On fluctuations in coin tossing, Proc. Natl. Acad. Sci. USA 35 (1949) 605--608.

\bibitem{Eu-Fu-Yeh-JCTA}
S.P. Eu, T.S. Fu, Y.N. Yeh, Refined Chung-Feller theorems for Dyck paths, J. Combin. Theorey Ser. A 112 (2005) 143--162.

\bibitem{Eu-Liu-Yeh-AAM}
S.P. Eu, S.C. Liu, Y.N. Yeh, Taylor expansions for Catalan and Motzkin numbers, Adv. Appl. Math. 29 (2002) 345--357.

\bibitem{Guo-WangXX}
V.J.W. Guo, X.X. Wang, A Chung-Feller theorem for lattice paths with respect to cyclically shifting boundaries. Preprint.

\bibitem{LiuSC-WangY-Yeh}
S.C. Liu, Y. Wang, Y.N. Yeh, Chung-Feller property in view of generating functions, Electron J. Combin. 2011, 18: P104.

\bibitem{Ma-Yeh-EJC}
J. Ma and Y.N. Yeh, Refinements of $(n,m)$-Dyck paths. European J. of Combin, 2011, 32: 92--99.

\bibitem{Narayana}
T.V. Narayana, Cyclic permutation of lattice paths and the Chung-Feller theorem, Skand. Aktuarietidskr. (1967) 23--30.

\bibitem{EC2}
R. P. Stanley, {\it Enumerative Combinatorics}, Vol. 2, Cambridge
University Press, Cambridge, UK, 1999.

\end{thebibliography}
\end{document}